# The Marquise du Châtelet: A Controversial Woman of Science

Dora E. Musielak
University of Texas at Arlington
dmusielak@uta.edu

**Résumé**. No woman of science has lived a more controversial life nor possessed a most contrasting character than Gabrielle Émilie Le Tonnelier, Marquise du Châtelet. On one hand, she was a woman of great intelligence, a philosopher of science, a student of mathematics, and she was an ardent supporter of Newton and his new laws of physics. At the same time, Émilie du Châtelet was an aristocrat society woman who gambled, enjoyed parties, and had several extramarital affairs, provoking numerous scandals in her native Paris. She was a passionate woman who was at ease conversing with the nobles at the court and with the most renowned scholars of her time. Émilie du Châtelet did not develop theorems nor discover new scientific principles. However, she studied mathematics with Maupertuis and Clairaut to better understand the geometrical language in Newton's *Principia*. In this article we review some important aspects of this controversial woman of science, exploring her relationship with the greatest scholars of her time.

**Résumé**. Ninguna mujer científica ha vivido una vida mas controvertida ni poseído un carácter mas contrastante que Gabrielle Émilie Le Tonnelier, Marquesa du Châtelet. Por un lado, era una mujer de gran inteligencia, una filosofa de ciencia, estudiante de matemáticas, y fue una ardiente partidaria de Newton y sus nuevas leyes de la física. Al mismo tiempo, Émilie du Châtelet era aristócrata, una mujer de alta sociedad que le gustaba apostar en juegos de azar, que disfrutaba de las fiestas, y tuvo varias relaciones extramaritales, provocando numerosos escándalos en su París natal. Ella era una mujer apasionada que conversaba fácilmente con los nobles de la corte y con los más reconocidos eruditos de su tiempo. Émilie du Châtelet no desarrolló teoremas ni descubrió nuevos principios científicos. Sin embargo, ella estudió matemáticas con Maupertuis y Clairaut para entender mejor el lenguaje geométrico en la *Principia* de Newton. En esta monografía tratamos unos aspectos importantes de esta mujer controvertida, explorando sus relaciones con los eruditos mas importantes de su tiempo.

## 1. Introduction

The Marquise du Châtelet, like other thinkers of the Enlightenment, aspired to understand the exact sciences, or at least to gain a reputation for understanding and contributing to expound the scientific principles derived from philosophical thought. She sought to understand the nascent science of mechanics through the philosophy of Newton, Descartes, Leibniz, and Wolff, and in the process she participated in the debate that raged in the eighteen century between Cartesians and Newtonians.

Madame du Châtelet never met Newton, as he died before she learned of his theories. However, she was perhaps one of the first women in history who tried to understand Newtonian physics. To contribute to the scientific enlightening of her contemporaries, du Châtelet made the first translation of Newton's *Principia* into French. Her intention was to make Newton's theories accessible and better understood, and in doing so, she contributed to the widespread campaign to promote Newtonian philosophy of natural science.

## 2. Her Early Life

Gabrielle Émilie was born to a wealthy family in Paris on 17 December 1706. Her father, Louis Nicolas Le Tonnelier de Breteuil, was a high-ranking official at the Court of Louis XIV holding the title of *Introducteur des Ambassadeurs*. He gave his daughter a superb education, hiring tutors to teach her Latin, foreign languages, horse riding, gymnastics, theater, dance, and music. Monsieur de Breteuil was a lover of literature, and so had a small literary circle that met periodically in his Parisian apartment. That is how Émilie met prominent intellectuals such as French author Bernard Le Bovier de Fontenelle and others interested in philosophy.

At age sixteen, Émilie was introduced to the Court of Versailles. That year, the young dauphin was crowned King of France Louis XV, known as *Louis le bien aimé*. When she was nineteen, Émilie married Florent-Claude, Marquis du Châtelet-Lomont, an army officer who was ten years her senior.[1] This was an arranged marriage, which became a burden for Émilie, for she had nothing in common with her husband. He was dull, while she was clever and witty. She enjoyed Paris society and life at court; he was more involved in his local regiment and had no intellectual interests. Of course her marriage to the marquis du Châtelet gave Émilie an even higher status in society, including *the droit du tabouret*, the honour of sitting (on a stool) in the presence of the queen,[2] and of travelling in her retinue. The queen was Marie Leszczyńska, daughter of the dethroned King of Poland, Stanisław Leszczyński. The royal wedding took place on 5 September 1725 (the same year Émilie married) at the Château de Fontainebleau. Undoubtedly, Émilie must have been part of the royal fête.

At the court of Versailles, Madame du Châtelet was a familiar guest, especially at the queen's card table. The queen enjoyed a game called cavagnole, a game of chance that Émilie played in excess, causing her to be in debt often. When the queen left for her autumn visit to Fontainebleau, du Châtelet had a place in one of the royal coaches. Her rank and precedence at other courts were undisputed.

Being a nobleman of hereditary rank whose illustrious family owned land in Lorraine, on the border of France, Émilie's husband was entrusted to defend and fortify against potentially hostile neighbors; he was away at war most of the time. Émilie and Florent-Claude had a daughter named Françoise-Gabriele-Pauline, born on 30 June 1726, and a son, Louis-Marie-Florent born on 20 November 1727. The following year, Émilie left their home in Semur-en-Auxois and moved to Paris, leaving her children behind entrusted to nannies, as was the custom of noble families. Madame du Châtelet orbited in the privileged circles surrounding the Queen of France,[3] and her friends came from the highest ranks of society. She attended the Court in Versailles and engaged in frivolous activities, including gambling and card playing, so popular among the aristocrats. Her social life was full of vanities, including the pursuit of love affairs.

---

[1] Hamel, Frank (1910). *An Eighteenth Century Marquise: a study of Émilie du Châtelet and her times*. London: Stanley Paul and Co. p. 25.

[2] Davidson, I., *Voltaire: A Life*, Pegasus, p. 102

[3] Polish-born Marie Leszczyńska (1703 – 1768) was the queen of France at that time. She was a daughter of King Stanisław Leszczyński of Poland (later Duke of Lorraine) and married King Louis XV of France in 1725. The longest-serving queen consort of France, Leszczyńska was grandmother of Louis XVI, Louis XVIII, and Charles X.

In 1732, the marquis du Châtelet returned from his military campaign and Émilie was summoned to their home in Semur where she became pregnant with their third child. When her husband rejoined his regiment to fight in the war of the Succession of Poland,[4] in January 1733, Émilie went back to Paris where her baby was born, named Victor-Espirit, on 11 April 1733. With that, she considered her marital responsibility completed, the baby was entrusted to nurses, and within weeks she resumed her carefree life full of whirlwind social amusements in Paris. Like many aristocrat women, Émilie du Châtelet was frivolous, fond of expensive clothes, wore diamonds and extravagant shoes, which must have cost a fortune to her husband.

## 3.    Madame du Châtelet and Voltaire

In May 1733, Émilie met the famous and controversial poet, philosopher, and playwright François-Marie Arouet, known by his *nom de plume* Voltaire. He was a *philosophe*, an intellectual of the Enlightenment famous for his wit and his attacks on the established Church. Because he was an outspoken critic of the government, Voltaire was frequently imprisoned or banished from Paris. When Émilie met him, Voltaire was 39 and had returned from exile in England where he lived from 1726 to 1729.

In England, Voltaire had become acquainted with the philosophical ideas of Bacon, Newton, and Locke. This experience inspired Voltaire to write *Lettres philosophiques* (translated as "Letters Concerning the English Nation"), twenty-five letters where he discussed his views on institutions, religion, philosophy, and the people in England. Voltaire praised the British parliament, which some considered to be an implied condemnation of the French monarch. He also discussed the new ideal of natural philosophy, the philosophical study of nature and the physical universe that was dominant before the development of modern science. In Letter 14, for example, Voltaire drew comparisons between Isaac Newton and French philosopher René Descartes. However, Voltaire was not a scientist and one can easily cast doubt on his understanding of Newtonian principles. In fact, a few years later Voltaire lost interest in any scientific research, or what he called "the sterile truths of physics," and went back to poetry and playwriting, creative areas in which he excelled.

In Paris, Émilie and Voltaire became inseparable, attending the opera, cabarets, theatres, and went together at the audiences of the Royal Court, forgetting the rules of decorum. Appreciating Émilie's interest in mathematics and natural philosophy, early in 1734 Voltaire introduced her to Pierre-Louis Moreau de Maupertuis, the most knowledgeable Newtonian scientist in Paris; Voltaire invited him to help Madame du Châtelet with her studies.

In 1734, Maupertuis began to give Émilie lessons in science and mathematics. Seduced from their first meeting, she became his mistress. She was overly passionate and pursued Maupertuis with the same ardor she studied. Even though she was already twenty-eight years old and married, Émilie behaved as an immature teenager and wrote many letters to Maupertuis declaring her affection for him. When her sixteen month-old child died, she wrote to Maupertuis begging him to come to console her. He did not reciprocate.

---

[4] The War of the Polish Succession (1733–1738) was a major European war for princes' possessions sparked by a Polish civil war over the succession to Augustus II, King of Poland that other European powers widened in pursuit of their own national interests. France and Spain, the two Bourbon powers, attempted to check the power of the Austrian Habsburgs in Western Europe, as did the Kingdom of Prussia; whilst Saxony and Russia mobilized to support the eventual Polish victor.

Although he remained a friend to her, Maupertuis became tired of this too pervasive love and asked Alexis-Claude Clairaut to continue teaching Émilie. Seven years younger than her, Clairaut was a prodigy mathematician known throughout Europe; he found in Madame du Châtelet an eager student and they became good friends. In 1736, Clairaut submitted to the Parisian Academy the lessons he gave to Châtelet, which were published with the title of *Élémens de Géometrie* (1753).[5,6] Year later, Clairaut also helped Émilie when she embarked in her research to understand Newton's theories.

In those years and throughout her life, Émilie wrote extensive letters to her friends and lovers, giving us a glimpse of her emotional states and her thought process. Those letters and the ones they wrote back reveal that she lived such a hectic social life (attending the opera, the theater, parties at the court), that it is difficult to grasp when she had time to study.

## 4.    Cartesian and Newtonian Ideas

In 1686, Isaac Newton published *Philosophiæ Naturalis Principia Mathematica*, introducing the laws of motion, and stating a new principle of gravitation to explain how the Solar System works. German mathematician Gottfried Wilhelm von Leibniz criticized Newton for not explaining how gravity works, while French scientists argued that the force of gravity was merely a supernatural idea. Many scholars had embraced René Descartes's theory of vortices, which postulated that space is entirely filled with matter in various states, whirling about the Sun. He also formulated a conservation of motion principle, asserting that the total quantity of motion—that is, the total quantity of mass times speed—always remains the same. Scholars who accepted Descartes's views were called "Cartesians."

In the 1730s, the scientific climate in Paris was infused with ideas from two opposite thinkers: Cartesians and Newtonians, i.e., those who followed Descartes theories and those who wanted to believe Newton. The Cartesians wanted to strip the Aristotelian universe of what they saw as its magic, its superstitious explanations, and substitute a mechanistic explanation in which we understood all problems of physics as problems of the communication of force by matter in motion to matter in motion. For the Cartesians, the Newtonian explanation of gravity as action at a distance —two masses with nothing between them affecting each other with gravitational pull— sounded like the Aristotelians and their secret, or occult, forces in nature. The ideas of Descartes and Newton were the topics of intense debate throughout Europe at that time.

Maupertuis was a Newtonian. During a six-month trip to England in 1728—a year after Newton died—Maupertuis had learned the theories in the *Principia* of Newton. In the autumn of 1729, Maupertuis went to Switzerland to study at the University of Basel where renowned mathematician Johann I Bernoulli welcomed him and became his mentor. Bernoulli had accepted Newton's results of universal gravitation but he also considered Leibniz's philosophical interpretations to provide an explanation for gravity, which Newton had left unexplained. Bernoulli taught Maupertuis calculus and introduced him to Leibniz's fundamental ideas of dynamics. At this time, Maupertuis became close friends with Johann II Bernoulli, Johann's I third son.

---

[5] Issu des leçons de Clairaut à la marquise du Châtelet, le manuscrit de C. 21 avait été présenté à l'Académie une première fois avant le voyage au Nord (cf. 18 avril 1736 (2)).
http://www.clairaut.com/n31aout1740po1pf.html
[6] http://www.e-rara.ch/zut/content/pageview/1330808

Back in Paris by July of 1730, Maupertuis began writing papers on mechanics, even if he had not fully understood Newton's inverse square law.[7] In 1731, Maupertuis was elected *pensionnaire géomètre* of the French Académie des Sciences The following year, he presented to the Academy memoirs on the laws of attraction, addressing the theories of Newton to demonstrate that gravitation is a universal physical principle. That year, Maupertuis published the treatise on the *figures des astres*[8] (Discourse on the various figures of the stars), supporting Newton's theory of gravitation, and opposing the ideas proposed by Descartes. When Johann II Bernoulli and his brother visited Paris in 1732, Maupertuis introduced them to the scientific circles, giving them access to a meeting of the Academy of Sciences the day when papers competing for a contest were read; one of the articles was written by their father, and the other by their brother, Daniel Bernoulli.

Having concluded his lessons to Émilie, Maupertuis returned to Basel in autumn 1734 accompanied by Clairaut. By now the relationship that Maupertuis had with the old Bernoulli was strained, but the friendship with his youngest son, Johann II, was blossoming. In Paris, the love affair between Émilie and Voltaire was on fire.

## 5.    Madame Du Châtelet and Voltaire at Cirey

When Voltaire's *Lettres* were published in France in 1735, the book was condemned, torn and burned. Once again, Voltaire had to leave Paris hurriedly, chased by the authorities. Madame du Châtelet came to his rescue and offered him a refuge in her husband's old Châteu de Cirey located 275 km east from Paris. The chateau was located near the border with Champagne and Lorraine, which was an independent province at the time, becoming an ideal sanctuary for Voltaire.

A few months later, Madame du Châtelet joined Voltaire. The freshly renovated Château de Cirey became the center of Émilie's scientific and philosophical work. She surrounded herself with books and scientific instruments of all kinds and was proud of its "pretty fine physics library, and telescopes." Madame du Châtelet and Voltaire had correspondence with many scholars, including Johann Bernoulli II, Maupertuis, Samuel König, Christian von Wolff, Charles Marie de La Condamine, and others who shared their intellectual interests.

Madame du Châtelet and Voltaire lived in Cirey for four years, dedicating much of their effort to studying and debating physical and metaphysical issues. She called Cirey "the land of philosophy and reason." It was evident to all that Voltaire and Émilie loved each other, away from the salons and the social whirl, away from the jealous glances and gossip. However, being both exuberant and sociable they invited many friends to Cirey; some guests stayed for weeks or even months.

---

[7] http://www-history.mcs.st-andrews.ac.uk/Biographies/Maupertuis.html
[8] Maupertuis, *Discours sur les différentes figures des astres, d'où l'on tire des conjectures sur les étoiles qui paroissent changer de grandeur, et sur l'anneau de Saturne, avec une Exposition abbrégée des systèmes de M. Descartes et de M. Newton*. Paris, 1732.

Visitors to Cirey wrote that Madame Châtelet and Voltaire gave lively parties in the evenings after engaging in their intellectual work. Madame de Graffigny, who visited the château between December 1738 and March 1739, described in her letters the theatrical frenzy of the inhabitants of Cirey. For one such event, "everyone was requisitioned to act. Even Madame du Châtelet's little daughter, Françoise Gabrielle Pauline, now twelve years old, was summoned from the Joinville convent, four leagues away. The little girl took naturally to Latin, and she learned her part in the play while travelling in the coach to Cirey."[9]

Italian writer Francesco Algarotti visited Cirey in 1735. During his stay, Algarotti wrote a book intended for the popularization of Newtonian philosophy addressed to women. With the title *Il Newtonianismo per le dame, ovvero Dialoghi sopra la luce*,[10] it was written as a series of dialogues. Algarotti wrote, "I am putting the last touches to my *Dialoghi* which have found grace in the eyes of the belle Émilie and the savant Voltaire. I try, when near them, to acquire those choice terms, that charming turn of speech with which I should like to embellish my work." Algarotti illustrated the book with an engraving of Émilie and himself set in a rustic scene, which represented the Cirey gardens with the chateau on the right. The marquise was, of course, highly flattered at being placed at the head of the work to represent "wit, grace, imagination, and science."

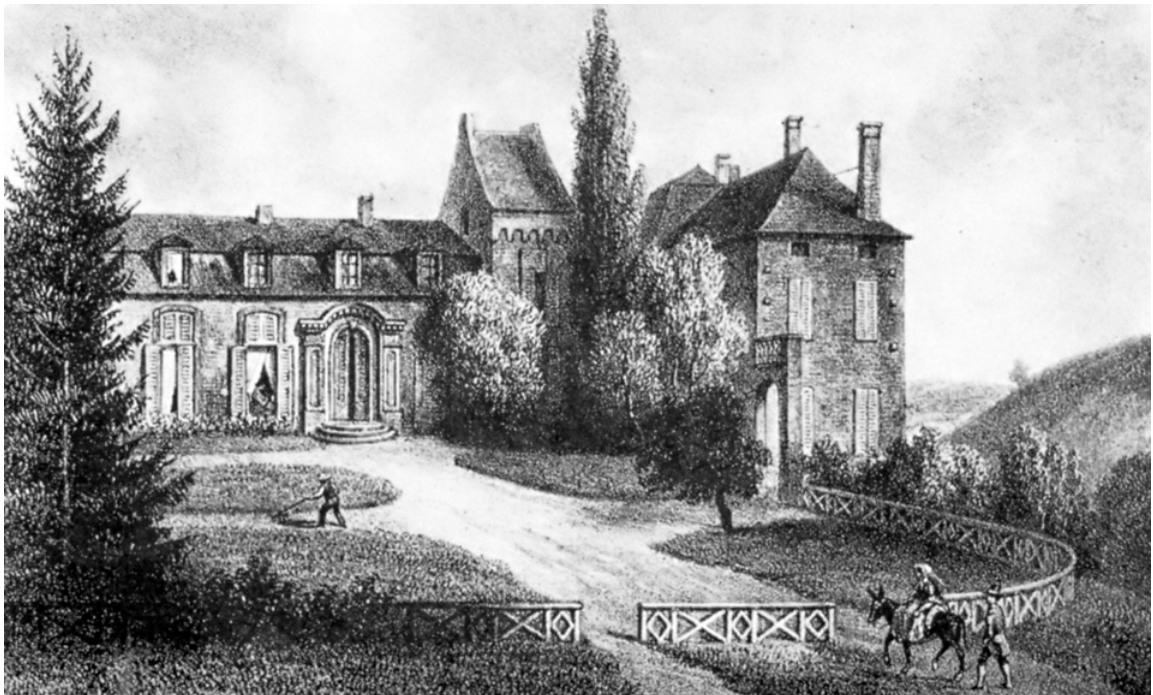

Image illustrative de l'article Château de Cirey. Lithographie d'époque
fr.wikipedia.org

---

[9] *Ibid*. p. 185. Based on the written description given by Mme. de Graffigny
[10] Newtonian [ideas] for the Ladies.

Despite her happiness at Cirey, Madame du Châtelet apparently also missed the gaiety of city life. In January 1736, she traveled to Paris where, "she had returned to her usual occupations, which were socializing, frivolity, philosophy, and — first of all in her heart — mathematics."[11]

Meanwhile, the dispute about the Newtonian theory of gravitation, and the proper methods for testing it, prompted Maupertuis to propose an expedition to Lapland, the polar region in Finland, to measure the oblateness of the Earth and prove what Newton predicted. In 1735, with the approval of King Louis XV, the Academy of Sciences sponsored Maupertuis to lead this scientific mission. He left France in 1736, accompanied by several young scientists, including Clairaut, Le Monnier, Camus and Celsius. On 13 November 1737, Maupertuis reported on his Lapland journey to the Academy, and the following year he published *Sur la figure de la terre*, the work which confirmed Newton's view that the Earth is an oblate spheroid slightly flattened at the poles. Maupertuis, now a very celebrated scholar, renewed his friendship with Madame du Châtelet and Voltaire and went to visit them in Cirey.

## 6.    Madame Du Châtelet on the Nature and Propagation of Fire

In 1737, the French Academy of Sciences issued a prize competition on the nature of fire and its propagation. Voltaire, perhaps overconfident from having read Newton and having learned some basic science from discussions with Maupertuis, Clairaut, and Châtelet, decided to submit a competing essay. He performed some experiments, but the philosopher could not achieve the required results. At the same time, Émilie was driven by her own scientific ambition, also prepared to submit her own entry for the Academy's contest, without telling Voltaire.

Madame du Châtelet wrote a 139-page essay that attempted to make a synthesis of all knowledge on the nature of fire. Of course, she didn't have new experimental facts nor did she offer a valid scientific theory. Her introduction is rather weak. She began her *Dissertation sur la nature et la propagation du feu* by saying that fire is difficult to define and it is impossible to fully understand its nature: "Le feu se manifeste à nous par des phénomènes si différents, qu'il est aussi difficile de le définir par ses effets, qu'il paraît impossible de connaître entièrement sa nature." She concluded the first part by saying that light and heat are two different and independent effects. To discuss the propagation of fire, she also addressed the nature of the Sun, saying "Le Soleil ne peut être un globe de feu," quoting Newton extensively to describe the colors of sunlight. Regarding the stars, for example, du Châtelet wrote:

"The light of the Sun appears to be yellow. Thus it is necessary that the Sun projects by its nature more yellow rays than others, ... Mr. Newton has shown in his Optics, page 216, that sunlight is abundant in this kind of ray.

---

[11] Hamel, Frank (1910). *An Eighteenth Century Marquise: a study of Émilie du Châtelet and her times*. London: Stanley Paul and Co. p. 68.

> It is very possible that in other systems, there are suns which project more rays of red, green, etc. and that the primary colors of suns that we never see are different from ours, and that, in short, there are in Nature other colors than those that we know of in our world."[12]

Some contemporary physicists use this passage to conclude that Madame du Châtelet had perceived a relationship between the position of a color in the spectrum and its energy, and foretold the existence of energies beyond the visible spectrum.

Be that as it may, Émilie did not discover the nature of fire or the laws that govern its propagation. The 1738 Academy's prize was awarded to 31-year-old Leonhard Euler,[13] and to two other lesser-known scholars.[14] With his *Dissertatio de igne* (the only entry written in Latin), Euler proposed that fire is the result of the bursting of tiny glassy balls of highly compressed air in the pores of bodies, so that "heat consists in a certain motion of the smallest particles of a body." Euler believed that all the phenomena associated with heat and fire could be deduced from the laws of mechanics without presupposing any "occult qualities." Euler stated, "light is the elastic vibration of the *ether* that is initiated by the explosions of little balls; hence, light is propagated by the same laws as sound."[15]

Du Châtelet and Voltaire were rather annoyed at the results of the competition. Was she expecting a favorable review since Maupertus was in the Academy's review committee (*commissaire pour les prix*)? They both complained, saying that it was unfair that their essays did not win. On 28 May 1738, Émilie wrote to Maupertuis: "Nous sommes au désespoir en voyant le jugement de l'Académie (des Sciences); il est dur que le prix ait été partagé, et que M. de Voltaire n'ait pas eu part au gâteau. Ce M. Fuller (sic) qui est nommé, est un Léibnitien, et, par conséquent, un Cartésien: il est fâcheux que l'esprit de parti ait autant de crédit en France."[16] Du Châtelet was referring to Leonhard Euler, who had won the highest accolades. What is surprising to me is that she seemed to refer to Euler with such disdain, and that she was naive enough to suggest that Voltaire had more merit than the young Euler! Clearly, Madame du Châtelet had no idea that Euler was the greatest mathematician of her time.

Moreover, Émilie du Châtelet was studying the philosophy of Christian Wolff and she knew that Euler was a critic of Wolffian ideas. Euler, the true scientist and most gifted mathematician, believed that knowledge is founded in part on the basis of precise quantitative laws, something that Leibniz's monadism[17] and Wolffian science were unable to provide.

---

[12] Pièces qui ont remporté le prix de l'académie royal des sciences de Paris, 1738. pp. 165-166

[13] Euler, L., *Dissertatio de igne, in qua eius natura et proprietates explicantur* (Dissertation on fire). Originally published in Pièces qui ont remporté le prix de l'académie royal des sciences de Paris, 1738, pp. 1-19. Opera Omnia: Series 3, Volume 10, pp. Republished in *Recueil des pièces qui ont remporté les prix de l'académie royale des sciences* 4, 1752, pp. 3-19 [E34a] 1738.

[14] The second winning essay was written by M. le Comte de Créquy, and the third by Pére Louis-Antoine Lozeran de Fiesc de la Campagnie de Jesus.

[15] Summary of Euler's ideas on the nature of fire from Clifford A. Truesdell's Introduction to Opera Omnia Series II, Volume 12. Copy of the essay is posted in Euler's Archive, http://eulerarchive.maa.org/

[16] *Lettres Inédites de Mme La Marquise du Châtelet, et Supplément a la Correspondance de Voltaire avec le Roi de Prusse, et avec différentes personnes célèbres*. Paris, 1818. p. 38.

[17] Monadism refers to an indivisible, impenetrable unit of substance viewed as the basic constituent element of physical reality in the metaphysics of Leibniz.

Du Châtelet was hurt that her own memoir did not win. Months later she wrote to Maupertuis explaining that despite the defects of her essay, she had hoped to attract the attention of the members of the review commission with her "boldness and novelty of her ideas."[18] Of course she didn't say that the Prize entailed a cash award[19] worth about as much as a modern Nobel Prize, approximately $ 1 million US.[20] Winning this prize would make the marquise quite wealthy and also enhance her scientific reputation.

Thanks to the favorable influence of their mutual friends, the essays of Madame du Châtelet and Voltaire were published in 1744 by the Academy of Sciences (together with the winning memoires), justifying it by stating that the "authors' names were likely to arouse the interest and curiosity of the public. One was by a lady of high rank, the other by one of the best of the poets."[21]

## 7.    Understanding and Explaining Newton

Newton's mechanics was a topic of study for Voltaire and Madame du Châtelet while they lived at Cirey. In 1738, Voltaire published his famous book *Elémens de la philosophie de Neuton* (Elements of Newton's Philosophy), which had the objective of providing a clear, accurate, and accessible account of Newton's philosophy. He published it as a *machine de guerre* directed against the Cartesians that, in his mind and that of his supporters, was holding France back from the modern light of scientific truth.  Some historians consider this the work that finally gave general acceptance of Newton's gravitational theories in France. With this work, Voltaire was viewed as an ardent activist of the new Newtonian physics. However, because of the breath of the topics exposed therein, it has been suggested that Émilie had much to do with the writing of this book. On 10 January 1738, she wrote to Maupertuis to thank him for his lessons, as this knowledge had allowed her to help Voltaire with his *Elémens*: "my companion in solitude has written an Introduction to Newton's philosophy which he dictates to me. I have the advantage over the greatest Philosophers in having had you for my teacher."[22] Voltaire admitted, *Minerva dictait et j'éceivais*.

 Not surprisingly, Voltaire dedicated the *Élémens* to the Marquise du Châtelet. He frequently praised her intelligence, saying she was a genius, calling her the Minerva of France, and a disciple of Newton, and he dedicated almost all of his work to her during their fifteen-year relationship. Du Châtelet responded in kind by writing a celebratory review of Voltaire's *Éléments* in the *Journal des savants*, the most authoritative French learned periodical of the day.

A second edition of the *Éléments* was published in 1745, and again later to include a new section devoted to Newton's metaphysics (Newton called it natural philosophy). At that time, scientific questions were addressed as a part of metaphysics. This is a division of philosophy that is concerned with the fundamental nature of reality that is outside of what can be perceived with the senses. Metaphysicians attempted to clarify the fundamental notions by which people understand the world, e.g., existence, objects and their properties, space and time, cause and effect, and possibility. Maupertuis, for example, was a mathematician who hoped to prove the existence of God by metaphysical principles and mathematics. These and other ideas he shared with Voltaire and Madame du Châtelet.

## 8.    Du Châtelet, Johann II Bernoulli, and Samuel König (Koenig)

On 12 January 1739, Maupertuis visited du Châtalet and Voltaire at Cirey and stayed for a few days. From Cirey he traveled to Basel to visit his friend, 29 year-old Johann II Bernoulli, and stayed in Basel until mid-March. He also met German mathematician Johann Samuel König.

König, who was twenty-seven at that time, had also studied with Johann I Bernoulli for three years, and also was taught by Daniel Bernoulli.[23] Shortly after returning from St. Petersburg in 1731, mathematician Jakob Hermann introduced König to the ideas of Leibniz. Impressed by Leibniz, in 1735 König traveled to Marburg to further his knowledge of philosophy and law under the guidance of Christian von Wolff. One of the most prominent German philosophers, Wolff wrote numerous works in philosophy, theology, psychology, botany, and physics. His series of essays all beginning under the title *Vernünftige Gedanken* ("Rational Ideas") covered many subjects and he especially expounded Leibniz's theories in popular form. Years later, König would become embroiled with Maupertuis and Madame du Châtelet in separate but equally bitter disputes about issues of priority and tutoring. But in 1739, associating with like-minded pupils and friends of the Bernoullis was an exciting prospect for both Émilie and Maupertuis.

Upon arrival at Basel, Maupertuis received a letter from du Châtelet saying that he  "inspired her to engage more closely with geometry and the differential calculus." She expressed frustration to having to learn alone, saying that she'd like to study under the guidance of one of the Bernoullis: Si vous pouviez déterminer un de MM. Bernoully à porter la lumière dans les ténèbres j'espère qu'il serait content de la docilité, de l'application et de la reconnaissance de son écolière. Je ne puis répondre que de cela, je sens avec douleur que je me donne autant de peine que si j'apprenais la calcul, et que je n'avance point, parce que je manque de guide.[24]

Maupertuis arranged for both his friend Johann II Bernoulli and Samuel König to accompany him back to Cirey in March. From her correspondence with Maupertuis, it seems that König spent the spring and summer of 1739 giving the Marquise du Châtelet lessons in mathematics and on the philosophy of Leibniz.  However, the relationship between the young tutor and the exuberant marquise soon became sour. Some have suggested that the rift was caused by a dispute over the payments for the lessons that

König gave Émilie.[25] According to Voltaire, she had been rather generous, providing him with an apartment in Paris, where König lived for several months. Certainly, if König had offered his services as a tutor he must have expected more than what she had paid him. In August, du Châtelet wrote to Maupertuis desperate to see him, mentioning that König was teasing her about it.

## 9.    Madame Du Châtelet's *Institutions de physique*

Having her essay on the nature of fire published by the Academy of Sciences gave Émilie a sense of accomplishment, and it helped to invigorate her scientific studies. She sought and obtained the approval to publish a treatise intended to present the philosophical principles of Newton. Her work was ready as early as 18 September 1738, but she delayed its publication, expecting to improve it by inclusion of the metaphysics of Leibniz. Finally in 1740, her *Institutions de physique* (Physical Institutions)[26] was published. This is a richly illustrated book on the principles of physics and mechanics based on Newton's work. It also expounds on the philosophical work of Leibniz and Wolff.

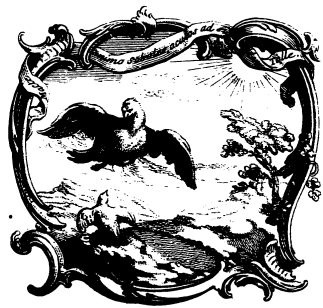

Frontispice of Émilie du Châtelet's *Physical Institutions* (1740)

The marquise du Châtelet claimed that the *Institutions* was intended for the education of her young son, Louis-Marie-Florent then thirteen, "whom I love with extreme tenderness." However, only researchers working on Newtonian mechanics could have appreciate her book. Du Châtelet covered the topics that were the subject of controversy at the time: the nature of forces. She gave credit to Leibniz for discovering the vital forces, i.e. to have "guessed one of the secrets of the creator."[27] Du Châtelet presented the *Physical Institutions* as a "small sketch" of the metaphysics of Leibniz and his disciple Wolff, stating: "I wanted to give an idea of the metaphysics of Leibniz that I confess to being the only one that has satisfied me, even though I still have doubts."[28]

Du Châtelet was proud of her work and sent copies of the *Institutions* to the most famous scholars of her time. She wanted to be recognized for her cleverness in the sciences as in her understanding of Newtonian philosophy. Driven by such desire, she also sent a copy of her work to Frederick II, King of Prussia, who was an intimate friend of Voltaire. Because Frederick did not like Émilie, he did not care for the gift. The king wrote to Voltaire to praise *his* work on Newtonian philosophy, but "When I read *Les Institutions de physique* of the marquise, I do not know if I have been deceived or if I am deceiving myself."[29]

Others didn't like Madame du Châtelet's book either or just couldn't appreciate her effort. In fact, the *Physical Institutions* was rather controversial. Du Châtelet had intended to write a synthesis of the systems of Descartes, Leibniz, and Newton.  In her opinion, Descartes had the merit of having handed over science to the forefront of knowledge and his method was revolutionary. Newton, meanwhile, demonstrated "using the most sublime geometry" that gravity is the cause of everything. In Leibniz's *The Theodicy* she found an answer to the question: how to reconcile the freedom of God and the necessity of the laws of nature? Leibniz believed that the entire universe comes from the aggregate of the monads, God, our soul, that of animals and particles of inert matter. The *monads* is Leibniz's idealistic view that all substances are simple un-extended substances. As exposited in his *Monadologie*, monads are elementary particles with blurred perceptions of one another. One can compare the monads with the corpuscles of the Mechanical Philosophy of René Descartes and others. In their view, the monads are the ultimate elements of the universe.

In Chapter XXI of the *Physical Institutions*, Madame du Châtelet questioned the theory of forces promoted by Jean-Jacques d'Ortous de Mairan. Mairan was an astronomer and physicist,[30] the Secretary of the Academy of Sciences at that time. Émilie assessed de Mairan's essay on estimation and measurement of the driving forces of the bodies (published in 1728), and criticized it because it advocated that the strength of the body corresponds to the simple product of mass velocity. Du Châtelet stated that, in determining the force of a body, the correct calculation was that of Bernoulli and Leibniz. In other words, the strength of the body should be the product of its mass by the square of its velocity. As soon as du Châtelet's book appeared in print in 1741, her critical assessment de Mairan caused a public dispute with him. De Mairan responded with a rather long letter regarding the question of the *forces vives*.[31]She fired back.[32]

---

[27] Letter from Du Chatelet to Maupertuis dated February 1738.

[28] Letter from Du Chatelet to Maupertuis dated 29 September 1738.

[29] Davidson, I., *Voltaire: A Life*.  Pegasus. p. 168.

[30] De Mairan's scientific work includes contributions to the theory of heat, observations of meteorological phenomena, and theoretical work on the orbital motion and rotation of the Moon. He noted a small nebulosity around a star closely north of the Orion Nebula in or before 1731, and mentioned in his *Traite physique et historique de l'Aurore Boréale* (Physical and Historic Tract of the Aurora Borealis), published 1733 in Paris and reprinted in the Journal des Scavans in 1754. Charles Messier later cataloged his nebula as M43.

At that time, the debate centered about the concept of *vis viva* (Latin for "living force"). The *vis viva* is an obsolete scientific theory that served as an elementary and limited early formulation of a conservation principle. It was the first description of what we now call kinetic energy but, as defined then, it was misunderstood and it caused a long and bitter debate regarding Newtonian mechanics. It had become clear to some scientists of the early eighteenth century that Newton's three laws were sufficient to explain the motions of "point masses." However, it was unclear how the laws could be applied to rigid bodies and fluids. The *vis viva* controversy started as a disagreement about Leibniz's concept and the ideas of Descartes. Leibniz's concept seemed to oppose Newton's theory of conservation of momentum, which Descartes advocated. Leibniz published in Acta Eruditorum a short note entitled "A Brief Demonstration of a Notable Error of Descartes and Others Concerning a Natural Law." In this note, Leibniz rejected the Cartesian equivalence between motive force, which Leibniz agreed is conserved in nature, and quantity of motion, which he argued is not. He had yet to mention *vis viva*.

Leibniz's note provoked heated exchanges with the Cartesians. The conservation of motion conceived by Descartes was difficult to abandon since they believed that all space is filled with matter. Leibniz agreed on the conservation of directional motion, but argued that because it is directional, unlike $mv^2$, and he conceded that it is not conserved in a collision. Then in 1695 Leibniz published a memoir on his new science of dynamics: the *Specimen dynamicum*. In this work, Leibniz introduced *vis viva* to differentiate between living and dead forces. As examples of dead forces he gave centrifugal force and gravitational or centripetal forces; Leibniz also included the forces involved in static equilibrium that, when unbalanced, initiate motion.

Referring to the metaphysical principle (that the effect must equal the cause), Leibniz calculated "the force through the effect produced in using itself up" and concluded that the force transferred from one equal body to another varies as the square of the velocity. For Leibniz the metaphysical principle established the priority of the conservation of living forces in changes of motion.

During the period 1676 − 1689, Leibniz concluded that in many mechanical systems (of several masses, $m_i$ each with velocity $v_i$) the quantity $\sum_i m_i v_i^2$ is conserved. He called this quantity the *vis viva* or living force of the system. This principle can be seen now to represent the conservation of kinetic energy in elastic collisions, which is independent of the conservation of momentum. However, many scholars of the eighteenth century did not understood this idea because they were influenced by Newton's prestige in England and of Descartes's authority in France, both of whom had promoted the conservation of momentum as a guiding principle. In other words, they believed the momentum $\sum_i m_i v_i$ was the conserved *vis viva*.

Thus, the issues that Madame du Châtelet grappled with required establishing the quantities that are universally conserved: Descartes's motion, Leibniz's *vis viva*, or what we now call momentum. The problem was controversial because it involved several other issues. One was semantic: the interpretation of "force" was not clear. The other was the metaphysical issue raised by Leibniz. And of course, the empirical issue of the apparent nonconservation of *vis viva* in the collision of soft bodies was even less clear. Nonetheless, because the concepts of momentum, work, and energy were not well understood, the idea of *vis viva* was somehow clouded by individual's understanding of mechanical principles. However, Émilie du Châtelet developed Leibniz's concept and, combining it with Gravesande's observations, showed that *vis viva* was dependent on the square of the velocities. Thus, one can say that *vis viva* is an early representation of kinetic energy. But it was not until early nineteenth century when *vis viva* started to be identified with energy.

---

[31] M. de Mairan, Secrétaire Perpétuel de l'Académie Royale des Sciences, a Madame ma Marquise du Chastellet sur la question des forces vives, en réponse aux objections qu'elle lui a fait sur ce sujet dans les *Institutions de Physique*.
[32] *Réponse de Mme *** à la lettre de M. de Mairan sur la question des forces vives*, J. de Savants, 1741.

It would take many years and new mathematical developments to resolve these issues. In fact, the study of motion carried out by Euler, and later by Lagrange, and others focused on specific problems and on principles from which their mathematical solutions could be derived.

## 10.    Accused of Plagiarism by Samuel König

After the *Physical Institutions* was published, Samuel König accused Émilie du Châtelet of plagiarism, suggesting that her ideas regarding Leibniz and Wolff were his own. It is true that König was very knowledgeable about Leibniz's natural philosophy and mathematics. When Maupertuis met him in Basel, Châtelet asked that König come to teach her mathematics. He came to Paris and during that interaction König initiated Émilie on the metaphysics of Leibniz.[33]

König charged that du Châtelet had copied sections of her *Institutions* from his own writings on Leibniz and Wolff. Before the accusations, the marquise herself had admitted as much. In the preface of the original *Institutions* (1740), in props XII, she wrote: "Je vous explique dans les premiers Chapitres les principals opinions de Monsieur de Leibniz sur la Métaphysique; je les ai puisées dans les Ouvrages du célèbre Wolf [sic], dont vous m'avez tant entendu parler avec un de ses Disciples, qui a été quelque temps chez moi, & qui m'en faisait quelquefois des extraits."[34] Of course she did not refer to König by name, and in subsequent editions she removed that statement, clearly denying that König had anything to do with her own assessment of the work of Leibniz and Wolff.

Du Châtelet and König had already a tense relationship that began when he was tutoring her. According to some accounts, trying to defend herself from the accusation she retorted that "he was hired to teach her algebra, not metaphysics." True, König replied, but since at the end of every lesson Émilie would remark smugly "Cela, c'est évident!" he was tired of her haughtiness and thus offered to "teach her truths of great importance but without one shred of evidence," meaning metaphysics. According to König, he succeeded so well in his teaching that her book included those lessons. She vehemently denied it, of course.

In a letter to Maupertuis on 21 Aug 1740, Émilie angrily demanded that he set the record straight: "On me mande de Berlin qu'il y passe pour constant que König me l'a dicté. Je n'exige sur ce bruit si injurieux d'autre service de votre amitié que de dire la vérité." What was Maupertuis supposed to do? Besides, in the same edition of the *Institutions*, the Note of the Publisher states that du Châtalet delayed its publication in 1738 in order to add the metaphysics of Leibniz.[35] Thus, even though she was brilliant and could form her own opinions, it would not be an exaggeration to say that du Châtelet in effect had learned from and took extracts from König's own writings to expand her own explanations.

---

[33] http://www-groups.dcs.st-and.ac.uk/~history/Biographies/Konig_Samuel.html
[34] Madame La Marquise du Châtelet, *Institutions de Physique*, Paris 1740. p. 12-13.

The *Institutions* presented the new science of motion in a manner that was accessible to non-scholars, providing her readers with the basic concepts needed to understand Newton' laws. However, there were no mathematical descriptions of these physical laws. Her approach was philosophical and her views were based on her study and on the observations of other scholars. Madame du Châtelet maintained a close relationship with the leading philosophers and had lively conversations about those topics in person and in her letters. In fact, she heavily quoted Clairaut. He also helped Émilie when she translated Newton's *Principia*, work where many of Clairaut's own theories appeared.[36]

The Marquise du Châtelet was unmistakably an intelligent woman. Being educated in the traditions of philosophy, science and geometry by Clairaut, Maupertuis, and König, she was a fine, forward thinker. However, her mathematical abilities cannot be compared with those of her mentors. Of course, her contribution to the popularization of Newtonian principles is very important. However, mathematical descriptions are essential to our understanding of physical phenomena. We require mathematics to providing a clear understanding of the relations which describe the interaction between forces and motion. That is why neither Madame du Châtelet nor her critics could fully address the phenomena in question, and neither one could develop new scientific theories. The debate was fruitless in that respect. At the time, however, there were no many people who were capable of well-founded judgment on the new science of mechanics. Therefore, du Châtelet must be forgiven if her lack of understanding on the purely analytical aspects of the problem led her to write her treatise as she did.

## 11.    Émilie du Châtelet and Euler

During her dispute with de Mairan, the marquise du Châtelet sought Euler's endorsement when neither Maupertuis nor Clairaut, her closest friends and mentors, came to her defense. Before that, Euler had written her a short congratulatory note dated 19 February 1740 (Julian calendar).[37] In this letter,[38] a very polite Euler, who was still in St. Petersburg, wrote to du Châtelet to thank her, although he had not received the *Pièce* that she had sent him through Maupertuis. I believed that the *pièce* referred to her *Institutions Physiques*, and if so then Maupertuis must have told Euler about it before it was published. However, du Châtelet herself did not know about it, because it was not until 1741 that she thought of sending Euler her work, hoping that he would help her appease her critics.

---

[35] *Ibid.* Avertissement du Libraire: Ce premier Tome des Institutions de Physique était prêt à être imprimé dès le 18 Septembre 1738 comme il paraît par l'Approbation, & l'Impression en fut même commencée dans ce temps-la ; mais l'Auteur ayant voulu y faire quelques changements, me la fit suspendre ; ces changements avoient pour objet la Métaphysique de M. de Leibnits, dont on trouvera une Exposition abrégée au commencement de ce Volume.

[36] http://www-history.mcs.st-and.ac.uk/Biographies/Clairaut.html

[37] Russia still used the Julian calendar while the rest of Europe used the Gregorian calendar. When it was 19 February in St. Petersburg, it was already March 1 in Paris.

[38] A photocopy of this note was provided by Rosanna Cretney, researcher at the Open University UK.

In any case, in 1740 Euler could not comment on her *Institutions* since he had not read it. However, he praised du Châtelet anyway, as he knew her work on the nature of fire published two years earlier (together with his dissertation) by the Parisian Academy of Sciences. Euler added: ... Il m'est bien glorieux, Madame, de me voir en lice avec une personne qui fait un des plus rares ornements de son Sexe, par le Lustre que Vous avez bien voulu répandre sur les sciences les plus relevées, en y portant la Sublimité de Votre génie. Cela seul, je le sens bien, est capable de me donner un Relief que je n'oserais espérer de mes faibles lumières ...[39]

In the summer of 1740, Frederick II ascended the throne in Prussia and immediately began efforts to hire Euler. On December 16, Frederick invaded Silesia, precipitating the War of the Austrian Succession. That same year, Euler won the Academy of Sciences Prize on the ebb and flow of the tides, sharing the award with Daniel Bernoulli and two other scholars. Frederick II summoned Maupertuis to Berlin to help him reestablish the Prussian Academy of Sciences. Maupertuis accepted and recommended Euler to lead the mathematical section of the Academy. At the same time, in Russia there was political unrest and open hostilities directed at the German community. These social and political conflicts may have contributed to Euler's acceptance of Frederick's invitation to join his Academy in Berlin.[40]

In the summer of 1741, the public debate between du Châtelet and de Mairan was raging. On 29 May, Émilie (still in Brussels) wrote to Maupertuis: Quelques intéressante que soit pour moi, Monsieur, ma dispute avec M. de Mairan, les nouvelles de ce qui vous touche m'intéressant bien davantage ; mais, comme la curiosité n'a nulle parte à cet intérêt, mon amitié est contente de vous savoir à Berlin en bonne santé. Cette lettre-ci vous prouvera combien vous m'avez causé d'inquiétude et de joie.

In June, after learning that Euler was on his way to Berlin, Madame du Châtelet wrote to Maupertuis that she wanted Euler to read her *Institutions*, and also her letter to de Mairan. In an undated and unfinished letter to du Châtelet,[41] Euler finally addressed her *Institutions*, concentrating in particular in her chapter on hypotheses. He was pleased to learn that she fought strongly and firmly with those philosophers who wanted to banish completely the hypotheses in physics that were, in his opinion, the only way to achieve a certain knowledge of the physical causes. Du Châtelet asserted in her book that physical laws had to be understood as hypotheses.

I believe Euler wrote this letter in response to Châtelet's inquiry when she was seeking corroboration of her ideas during her dispute with de Mairan in 1741. In fact, on August 8, Émilie wrote to Maupertuis (she was still in Brussels), after learning from the newspapers that Euler had arrived to Berlin. She wrote: "Les gazettes disent Euler à Berlin; cela est-il vrai? Est-ce vous que l'y avez attiré; je ne sais s'il ne s'en repentira pas: il est vrai qu'il vient de Petersbourg; mais il y a bien des façons de perdre au change. Je voudrais lui envoyer les Institutions et les pièces de ma dispute avec Mairan: pourriez-vous lui faire tenir, je vous les ferais remettre ?" Euler had arrived in Berlin at the end of July 1741.

Euler responded to Châtelet as follows: "Mais j'espère qu'une bonne partie de ces gens changeront bientôt leur sentiment après avoir lu Votre admirable dissertation sur les hypothèses: et je ne doute nullement, que Mr. Demairan ne soit entièrement convaincu par les solides raisons que Vous avez opposées à ses idées si mal fondées sur la force des corps. Au reste je plains fort, que cette matière ait été déjà si longtemps le sujet d'une dispute si forte, par laquelle la mathématique a perdu beaucoup de sa réputation."[42]

Then Euler shared his thoughts regarding the concept of force, starting with the first law of Newtonian mechanics. In Euler's own words: "Je commence par le premier principe de la Mécanique que tout corps par lui même demeure dans son état ou de repos ou de mouvement. A cette propriété on peut bien donner le nom de force, quand on ne dit pas que toute force est une tendance de changer l'état, comme fait Mr. Wolf. Touts corps est donc pourvu ... "[43] Unfortunately, the letter ends at this point. However, the incomplete sentence suggests that Euler wanted to expand on the nature of forces.

Let's remember, Euler was an expert in mechanics, and his mathematical prowess was firmly established. Euler was the first to publish in 1736 a monumental book in which he applied the full analytic power of the calculus to Newtonian dynamics. Although Newton's *Principia* was fundamental to the understanding of mechanics, it lacked in analytical sophistication, i.e., the mathematics required to explain the physics of motion was rather obscure, as Newton had preferred to use geometrical arguments rather than his own invented calculus. Euler, on the other hand, rewrote the *Principia* using analysis so that the physics could be better understood. With his two-volume *Mechanica*, Euler described analytically the laws governing movement and presented Newtonian dynamics in the form of mathematical analysis for the first time. Later, he produced another mechanics text in which he introduced the equations for rigid body rotations.

We do not know what Châtelet responded to Euler. Had she met him in person, it is doubtful that Émilie would have liked Euler. For one thing, Euler was a pious, simple family man, too bourgeois, too unsophisticated and unpolished for her taste. Émilie du Châtelet preferred men like Voltaire —literary, witty, sophisticated—or elegant scholars like Maupertuis and the young Bernoulli, men who were charmed by her own wittiness. Euler and Madame du Châtelet belonged to two completely different worlds.

## 12.    End of a Love Affair and a New Beginning in Newtonian Mechanics

At age thirty-seven, the state of Châtalet's personal life was in shambles. In the past few years, her love affair with Voltaire had entered a period of stress and alienation. Although they continued to live (most of the time) and travel together, their feelings and interests were pulling them in different directions. Voltaire went away, for months at a time, and she wrote him letters full of reproaches. Even when they were together, the couple quarreled intensely and often.

---

[42] *Ibid.*
[43] *Ibid.*

In 1743, Madame du Châtelet arranged for the marriage of her daughter Pauline to the old Duke de Montenero-Caraff; the seventeen-year old girl went from the convent to her husband's residence in Capodimonte in Italy and never saw her mother again. In the autumn of the same year, Voltaire visited Frederick in Berlin, staying for several weeks. The King of Prussia did not want to host Madame du Châtelet so she stayed behind; Voltaire's letters from Berlin grew fewer and less loving. By then she must have began to get alarmed at the coolness of her lover.

The love relationship between Émilie and Voltaire finally ended when she discovered that Voltaire was having an affair with another woman. Surprisingly, they continued living and travelling together. At the same time, her former lover, Maupertuis, was in Berlin where he arranged his own marriage to Eleonor Borck. In the spring of 1745, Maupertuis returned to Paris to tidy up his affairs and announce his wedding, which took place on 25 August 1745. This must have hurt and humiliated Émilie, who had no other love prospects. Du Châtelet needed to seek an outlet for her loneliness and for the rejection by Voltaire and Maupertuis. She returned to her studies of Newton's philosophy. Her friend Clairaut was there, willing to help with her next project.

## 13. Du Châtelet and Her Translation of Newton's *Principia*

In 1745, Émilie du Châtelet began the greatest work of her life: the translation into French of the Latin text of Newton's *Principia*. She desired to be recognized by her own scientific skills, which were subject to mockery by many. At the same time, she wanted to put an end to the Cartesian's ideas, which were still embraced by some scholars.

On 20 December 1745, Clairaut approved the order for the translation, and on February and March of 1746, the royal decree document was signed by representatives from the king, from the Royal Chamber of Booksellers and Publishers of Paris, and by Breteuil du Châtellet. Émilie was ready to carry out a project that would make her name known throughout Europe.

The Latin of Newton must have been quite difficult to translate and it required that Newton's assertions be verified and explained further. Clairaut met with Émilie frequently to help her with the required calculations.

At the end of the first volume, which covered Newton's Books I and II, du Châtelet included a *Table Alphabétique*, a sort of Index in which she defined the terms used and referred them to particular sections of the book. For example, she defined the term "gravity" as follows:

*Gravite, elle est d'un autre genre que la force magnétique. III. 21 Cor. 3. Elle est mutuelle entre la terre & ses parties. III. 32. à la fin.*
*Sa cause n'est point assignée. III. 178 à la fin.*
*Elle a lieu dans toutes les planètes. III. v. Cor. I.*
*Et elle décroît hors de leur superficie en raison doublée de la distance au centre. III. viii.*
*Et de leur superficie vers le centre, elle décroît dans la raison simple des distances à peu près. III. ix.*
*Elle a lieu dans tous les corps, & elle est proportionnelle dans chacun d'eux à leur quantité de matière. III. vii.*

*C'est par la force de la gravité que la Lune est retenue dans son orbite. III. iv.*

The second volume was devoted to Newton's *System of the World*. As Newton himself stated in the introduction (and du Châtelet finely translated it), in the preceding books he had laid down the principles of philosophy; principles not philosophical but mathematical ... these principles are the laws and conditions of certain motions, and powers of forces, which chiefly have respect to philosophy .... Newton illustrated these principles in the System of the World, starting with four "rules of reasoning in philosophy," followed by descriptions of phenomena related to the motion of planets.

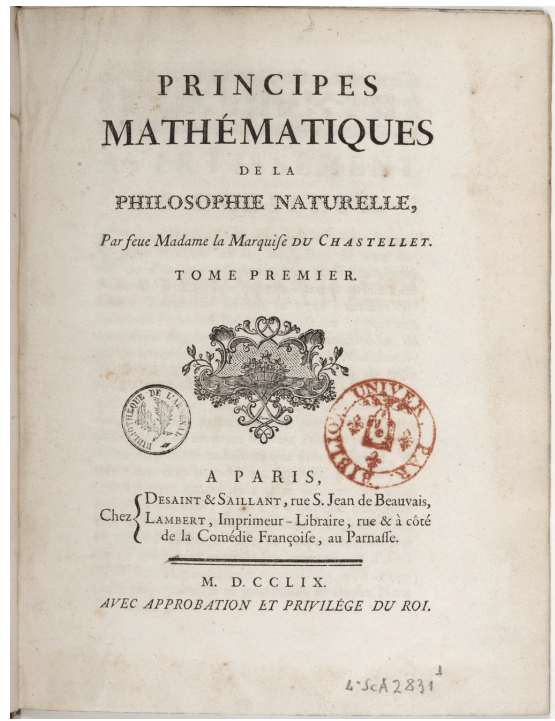

Frontispiece of Émilie du Châtelet's Translation of Newton's *Principia* (1759).

The translation of Newton's Book III ended with the General Scholium (pp. 174-180). Then du Châtelet followed with her explanations of Newton's work: " Short exposure of the System of the World: Explanation of the main astronomical phenomena from the principles of Mr. Newton. [44] Her commentary began with a brief historical overview of the discoveries of astronomy that culminated in the Copernican model and Kepler's laws of planetary motion, which led Newton to derive the general laws of motion and gravity. She described the solar system, defined the terms used, offered quotes from Newton and other scientists. This commentary was followed by a more scholarly part titled "Analytic solution of the major problems concerning the system of the world," which has the trademark of Clairault.

---

[44] *Exposition Abrégée du Système Du Monde: Explication des Principaux. Phénomènes astronomiques tirée des Principes de M. Newton.*

Du Châtelet considered the shape of the orbits, planets, in the different hypotheses of gravity, of light, of the figure of the Earth. For example, while discussing the determination of the shape of the Earth, Émilie quoted Clairaut's own work noting that Newton assumed that the Earth was an elliptic spheroid with homogeneous mass and that Clairaut "had checked the first of the assumptions was correct and that Newton had neglected,"[45] even though it is very important to ensure the true relation of the Earth's axes. In 1743, Clairaut had published *Théorie de la figure de la Terre*, confirming the Newton-Huygens belief that the Earth was flattened at the poles. From this point on, du Châtelet quoted extensively from Clairaut's work.

The last volume ended with a summary of the work of Daniel Bernoulli on tides. In 1740, Bernoulli, Euler and Colin Maclaurin[46] had won jointly the Paris Académie des Sciences prize for a study of the tides. However, Madame du Châtelet liked Bernoulli's work better, finding it clearer because she wrote on p. 261: *Je me fuis sur tout attachée à lire celle de M. Bernoulli, dans laquelle il m'a paru trouver plus d'ordre, de netteté & de précision ; & j'espère que le Lecteur me saura gré, si pour Commentaire de notre Auteur sur cette matière, je lui donne un abrégé du Traité de M. Bernoulli.* And thus she took Daniel Bernoulli's Treatise on Newton's Theory of Tides as the basis for this section.

## 14.    The End of Madame du Châtelet

Late in 1748, Madame du Châtelet and Voltaire traveled to Lunéville to visit Stanisław I Leszczyński, former King of Poland and father of the Queen of France. Since 1737, Leszczyński was the Duke of Lorraine, living in his château in Lunéville, located about 365 km from Paris. The grand Château de Lunéville, built in the style of Versailles in 1702, was the residence of the duke of Lorraine. A devout catholic, author, and philanthropist, Stanisław Leszczyński had a church built and several follies in his gardens for the amusement and education of visiting Polish nobility and followers of the Enlightenment, including famous philosophers and men of letters. Leszczyński's court was a place of frivolity and amusement where idle people dominated social life.

At the Château de Lunéville, Madame du Châtelet met Jean François, marquis de Saint-Lambert, with whom she fell madly in love. Émilie was forty-two years old, and the man, age thirty-one, was a mediocre poet, uninteresting and quite conceited. Almost immediately, Émilie began flirting with Saint-Lambert, and became infatuated. She showered him with an increasing stream of passionate letters, exchanged across her boudoir and his bedroom at the château.

At first Voltaire pretended not to notice, but soon enough the love affair was obvious to all around. It has been written that during this time, du Châtelet lived the last great passion of her life.[47] She had had many lovers, but she succumbed to a mad passion

---

[45] « En déterminant le rapport des axes de la terre, M. Newton outre la gravité mutuelle des parties de la matière, a encore supposé que la terre était un sphéroïde elliptique, & de plus que la matière était homogène. M. Clairaut, dans son Livre de la figure de la terre, a fait voir que la première supposition était légitime, ce que M. Newton avoir négligé de faire, quoique cela soit fort important pour s'assurer qu'on a le vrai rapport des axes de la terre. » p. 61.

[46] Scottish mathematican Maclaurin (1698-1746) published the first systematic exposition of Newton's methods, written as a reply to Berkeley's attack on the calculus for its lack of rigorous foundations.

[47] For an account of the love affair, see D. W. Smith, *Nouveaux regards sur la brève rencontre entre Mme Du Châtelet et Saint-Lambert*. In The Enterprise of Enlightenment. A Tribute to David Williams from his

for Saint-Lambert, a man whose intellectual interests were far beneath hers. Despite reproaches, disputes, marital strife, the love of Émilie for Voltaire had been deep, passionate, based on the sharing of knowledge. For fifteen years, they were together. However, during her time with Saint-Lambert, she lost any insight, any autonomy. For Saint-Lambert, Émilie sacrificed her position, her future, her children, her husband, and her best friend.[48] She put aside her translation of Newton while she wooed her new love interest.

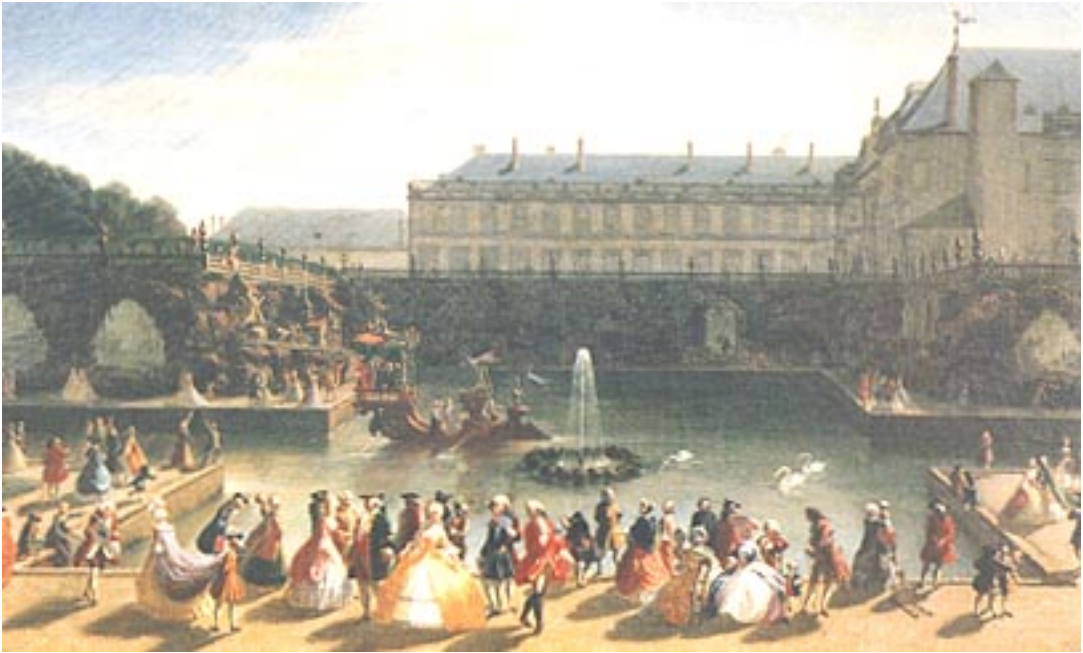

Château de Lunéville during one of Stanisław Leszczyński's fêtes.

Du Châtelet and Saint-Lambert spent three months together, until Voltaire forced her to return to Paris. By then it was too late; Émilie was pregnant. Hamel wrote that Émilie was "a woman to whom to love and to be loved is of the utmost importance, and who is prepared to sacrifice all else in life, honour included."[49] Madame du Châtelet first told Voltaire about the pregnancy and then to Saint-Lambert. With the former she pondered how she would give the news to her husband. Florent-Claude du Châtelet was summoned to Cirey where the two lovers (Voltaire and Saint-Lambert) came to be part of the revelation. What kind of man was her husband who seemingly accepted such a scandalous affair?

While pregnant Madame du Châtelet continued working on her translation of Newton's *Principia*. In Paris, Clairaut visited her every day to help her. In June 1749, Émilie terminated her work with Clairaut but the book was yet to be completed. She and

Voltaire left Paris, because she wanted to have her baby at the château in Lunéville. She arranged to have the queen's suite. That's how much influence she had on the Duke of Lorraine.

The trip from Paris to Lunéville must have been especially difficult for the six-month pregnant Émilie. Travelling in the eighteenth century was rather difficult and strenuous. Even if she travelled in some comfort —in an upholstered carriage accompanied by footmen and liveried coachmen and maids—the 365 km voyage from Paris to Lunéville required several days. The speed of travel depended on the condition of the roads, the terrain, the season and the weather. In good weather, the unwieldy vehicle, drawn by four horses, would lurch along the rutted roads at an average speed of about 5.5 to 7.5 km per hour, and it would need to stop for the exchange of horses after 3.5 or 4.5 hours. The exchange in the post stations often took up to two hours. After a rain, carriages often became stuck in the mud. Moreover, the accommodations along the way were often poor, with humid rooms and insufficient hygiene being quite normal.

On 4 September 1749, Émilie gave birth to a baby girl after an easy delivery attended by the King's doctor. The baby was baptized Stanislas-Adélaïde du Châtelet. After the birth, the mother seemed to be perfectly well, and the baby was given to nurses. Then Émilie contracted a fever. On 10 September two more doctors were called from Nancy, who tried in vain to restore her health. Émilie du Châtelet died that night, surrounded by her husband, Florent-Claude du Châtelet, and by her two lovers, Voltaire, and Saint-Lambert.

According to Voltaire, the baby herself died a few days later, regretted by nobody.[50] Had she lived, Émilie's daughter would have been shut up in a convent and never heard of again. Unless she would have inherited her mother's intellect and her indomitable free spirit.

## 15.  Émilie, the Woman

**Appearance**. People who knew her described Madame du Châtelet as having a striking appearance, rather than beautiful, lively rather than lovely. She was tall, had sea-green eyes, and a dark complexion. Émilie was feminine, enjoyed dressing up in lavish gowns adorned with pompons in the latest styles; she loved showing off in public glittering with diamonds. Contemporary portraits show that Mme du Châtelet was a pretty, good-looking woman. And no doubt, due to her intelligence, aristocratic breeding, her fine intellect and superb education, Madame du Châtelet must have appeared stunning and rather fascinating to those who conversed with her.

**Personality**. From hers and Voltaire's extensive correspondence, and also from various testimonials of their contemporaries, one can form a rather accurate image of Émilie du Châtelet's personality.  Of course, Voltaire's assessment of *la divine Émilie* was tinted with the rose coloured glasses of his love for her. The marquise du Châtelet was a complex, uninhibited woman. Her character can only be described as an extraordinary contrast between her intellectual powers and her enthusiasm for engaging in trivial, frivolous activities. Émilie was an enthusiastic singer of opera[51] and she loved acting. But by all accounts, she comes across as brilliant, deeply emotional, sophisticated, sociable, and a rather snobbish person.

---

[50] Some sources give the date of the infant's death as 6 May 1751 in Lunéville.

Émilie du Châtelet was arrogant, pretentious, vain, and self-assured. Feelings of superiority were rather innate in her. These personality traits were to help her deal with the criticisms of her work. Feelings of resentment, due to criticism or being overlooked, were familiar to du Châtelet. When de Mairan disputed her assertions in the *Institutions*, and others criticized it, she fought hard and made public her own views, resenting those who did not come to her defense. König accused her of plagiarism, for example, but she firmly denied it.

Her famous quote sums up what Émilie du Châtelet thought of herself and others: "Judge for my own merits, or lack of them, but do not look upon me as a mere appendage to this great renowned scholar or that, that I am in my own right a whole person, responsible to myself alone for all that I am, all that I say, all that I do. It may be that there are metaphysicians and philosophers whose learning is greater than mine, although I have not met them. Yet, they are but frail humans, too, and have their own faults; so, when I add the sum total of my graces, I confess that I am inferior to no one." She wrote this passage in a letter to Frederick of Prussia, a man who openly disliked her.[52]

Despite possessing a fine intellect, Madame du Châtelet lived a reckless life. She loved going to parties and to the theater, to dances and to the Opera. Émilie du Châtelet was known for her financial irresponsibility and her heavy gambling. By the end of 1743, she was in debt, penniless, and Voltaire was tired of bailing her out. In January 1744, Émilie wrote asking him for financial help (a loan for 50 louis): "I have sent £500 to M. du Châtelet for kitting out his son.[53] I shall pay you back in rental for the house, or else, if you prefer, here is the note from M. du Châtelet which, luckily, I did not tear up ... keep it and lend me the money ..." Voltaire was angry and told his acquaintances he was tired of her dissolute social life. When she died, she left a great number of debts, which her husband must have paid.

Émilie du Châtelet knew how to exert her power to gain favors. She lobbied for the military promotions of her husband and son. She interceded with the authorities on Voltaire's behalf and wrote glowing commentaries to promote his work.

Her maternal instincts were not developed. It has been written that she acted according to her aristocratic upbringing and the mores of her day, and she was rather influenced by the ideas of the French Enlightenment. Only once did she express sorrow for her children when her third baby died. However, because she wrote that letter to Maupertuis during the time she was ardently pursuing him, it is possible that the expression of grief in that midnight letter was a ploy to get him back.

Aside from that, one wonders about du Châtelet's thoughts regarding the scientific education of women. Du Châtelet dedicated her *Institutions* to her young son: "You are, my dear son, in that happy age at which the mind begins to think, and in which the heart has as yet none of the keener passions to trouble it. It is perhaps the only time in

---

your life that you could give to the study of nature; soon the passions and pleasures of your age will carry all your moments and when this passion of youth will be passed, and that you have paid to the intoxication of worldwide the toll of your age and your state, the ambition will take hold of your soul; and even in this age and that often is not more mature, you would like to apply to the study of real science, your mind having more flexibility which is the sharing of the beautiful years, you would need to buy a painful study what you can learn today with extreme ease."

Was the Marquise du Châtelet prejudiced regarding the education of girls? She sought for herself and her son a science education, but not for her daughter. Why? Émilie knew about Italian scholars Laura Bassi and Maria Gaetana Agnesi. In fact, on 1 April 1746, she was elected and inscribed in the register of members of the Italian Academy of Sciences in Bologna. Very proud of this recognition, du Châtelet was with Bassi and Gaetana part of a tiny group of women scholars in the eighteenth century.

Yet, Émilie readily accepted that, for well-breed girls, minimal convent schooling was sufficient. In fact, Madame du Châtelet sent her daughter to a convent to receive the basic education that girls of the nobility typically received before marriage. She wrote: "I have always believed that the most sacred duty of men was to give their children an education which hinders them in after years from regretting their youth, which is the only time in which they can really learn."[54] Clearly, she was referring to male children, and it implies that Émilie did not consider that her fourteen-year-old daughter required a scientific education. However, the son did not follow on the footsteps of his scholarly mother; Louis-Marie-Florent sought a career in the military, just as his father; at sixteen, the young man joined the army as an officer cadet.

## 16.    Les Passions d'Émilie

Émilie du Châtelet loved men, two or three at the same time, and so easily she transferred her affections among them. Her first lovers included the Comte de Guébriant, for whom she attempted suicide after he left her. In 1729, she had an affair with Armand de Vignerot du Plessis, the Duke of Richelieu, a womanizer aristocrat, and although their relationship was short-lived, they continued to be frequent correspondents for over a decade. Émilie also had an affair with Maupertuis, her teacher of mathematics who tried in vain to keep a distance from her. Then in 1733 she met and became the most famous lover of French philosopher Voltaire. He admired her dearly for her intellectual passion and called her "*la divine Émilie*," dedicating many of his works to her.

Voltaire and Madame du Châtelet lived together for many years, in her husband's estate house. What did the marquis du Châtelet think of such arrangement? He seems to have accepted Émilie's independence and social freedoms. For some it may be difficult to believe that the husband did not object to her scandalous affairs with Voltaire and Saint-Lambert. Eighteenth-century manners may have been free and easy, but du Châtelet's behavior and her husband's acceptance are difficult to understand, even today.

Voltaire was deeply affected by the death of his divine Émilie. On 15 October 1749, he wrote to Frédéric II, King of Prussia: "I lost a friend of 25 years, a great man

---

who had the defect to be a woman, and that all of Paris regrets and honors." Saint-Lambert soon recovered from his grief and moved to Paris around 1750.

Émilie du Châtelet was an emotional, passionate, an exuberant woman. It is difficult for one to fully grasp her fiery and passionate temperament and her desperate craving for men's love. She lived a reckless social life, she was devoid of maternal feelings, and was rather inconsiderate when it came to her husband and her children. Madame du Châtelet did not express sadness for not seeing her daughter after marrying her to Alfone Caraffa Despina, duc de Montenero-Caraffa in 1743. Françoise-Gabrielle Pauline went to live with her husband and never returned to Cirey or Paris. She died in 1754. Émilie's son, Louis-Marie-Florent became an army officer, and later ambassador to Vienna and London. He succeeded his father as marquis du Châtelet, and later became a duke. He was executed during the French Revolution in 1794; he was sixty-seven.

## Postscript

In 1759, ten years after her death, Clairaut published Émilie du Châtelet's work as *Principes mathématiques de la philosophie naturelle, par M. Newton* (traduit du latin) par feu Madame la marquise du Châtelet.[55] Voltaire wrote the preface. The book included two volumes in- 4°: the translation of Newton's *Principia* (1726 third edition) filled a volume and a half. The last half of the volume (about 100 pages) contains a "short exposition of the system of the world and the explanation of main astronomical phenomena from the principles of Newton." The book ended with a Glossary Index, where Émilie expanded on many terms that were new or important such as definitions of centripetal and centrifugal forces, and gravity. This first translation of Newton's *Principia* became a popular reference in France.

Although Émilie du Châtelet's efforts were limited to commentary and synthesis, her publications contributed to the dissemination and learning of Newtonian science in the middle of the eighteenth century.

The Marquise du Châtelet was buried at the Église St. Jacques in Lunéville. Today, we find a humble marble tombstone, an almost imperceptible reminder that her remains are there. Upon entering the church, one can almost stumble over the simple dark grey stone lying at the base of a vestibule column bearing the simple inscription: Sous cette dalle de marbre noir repose (under this black marble slab rests) Gabrielle Émilie Le Tonnelier de Breteuil Marquise du Châtelet, Paris 17 décembre 1706, Lunéville 10 septembre 1749 – Femme de Sciences et Philosophe.

## Acknowledgement

My sincere thanks to Rosanna Cretney, researcher at the Open University UK, for sharing with me the photocopy of Euler's letters to the Marquise du Châtelet. I am also grateful to Professor Dominic Klyve, Central Washington University, for his help in locating the copy and for his enthusiastic encouragement of this work.

---

[55] http://gallica.bnf.fr/ark:/12148/bpt6k1040149v